\documentclass[a4paper,12pt]{extarticle}
\usepackage{amsmath, amsthm, amssymb}

\usepackage{lscape}
\usepackage{fullpage}
\usepackage{csquotes}
\usepackage{cite}

\usepackage{algorithm, algorithmic}
\usepackage{epstopdf}
\usepackage{epsfig}
\usepackage{dsfont}
\usepackage{enumitem}
\usepackage{flexisym}
\usepackage{amsmath}
\usepackage{color}
\usepackage{algorithm, algorithmic}
\usepackage[section]{placeins}
\usepackage[lofdepth,lotdepth]{subfig}
\usepackage{mathtools}

\DeclarePairedDelimiterX\Basics[1](){ #1}
\usepackage[mathscr]{euscript} %euler for tensors
\usepackage{mathtools}
\providecommand{\keywords}[1]{\textbf{\textit{Keywords and phrases}} #1}

\DeclarePairedDelimiter{\floor}{\lfloor}{\rfloor}
\usepackage{blindtext} % for dummy text

\theoremstyle{plain}

\newtheorem{thm}{Theorem}
\newtheorem{lem}[thm]{Lemma}
\newtheorem{cor}[thm]{Corollary}
\theoremstyle{definition}
\newtheorem{dfn}[thm]{Definition}

\theoremstyle{remark}
\newtheorem{rem}[thm]{Remark}

\newcommand\onet{\mathbb{U}}
\newcommand{\m}{\text{max}}
\usepackage[toc,page]{appendix}
\usepackage{indentfirst}

\title{%\vspace{-1.5cm}            % Another way to do
Learning tensors from partial binary measurements}
\author{Navid Ghadermarzy, Yaniv Plan, \"{O}zg\"{u}r Y{\i}lmaz\\
Department of Mathematics\\
University of British Columbia}
\begin{document}
\maketitle
%  \tableofcontents\newpage
\begin{abstract}
In this paper we generalize the 1-bit matrix completion problem to higher order tensors. We prove that when $r=O(1)$ a bounded rank-$r$, order-$d$ tensor $T$ in $\mathbb{R}^{N} \times \mathbb{R}^{N}  \times \cdots \times \mathbb{R}^{N}$ can be estimated efficiently by only $m=O(Nd)$ binary measurements by regularizing its max-qnorm and M-norm as surrogates for its rank. We prove that similar to the matrix case, i.e., when $d=2$, the sample complexity of recovering a low-rank tensor from 1-bit measurements of a subset of its entries is the same as recovering it from unquantized measurements. Moreover, we show the advantage of using 1-bit tensor completion over matricization both theoretically and numerically. Specifically, we show how the 1-bit measurement model can be used for context-aware recommender systems.
\end{abstract}
\section{Introduction}\label{introduction}
The problem of \textit{matrix completion}, i.e., recovering a matrix from partial noisy measurements of a subset of its entries, arises in a wide variety of practical applications including collaborative filtering \cite{goldberg1992using}, sensor localization \cite{singer2010uniqueness,biswas2006semidefinite}, system identification \cite{liu2009interior} and seismic data interpolation \cite{aravkin2014fast}. The low-rank structure of a matrix makes it possible to complete it by sampling a small number of its entries (much smaller than its ambient dimension). Matrix completion is useful for applications where acquiring the full data is either expensive or impossible due to physical limitations.\newline

Assuming that we have access to $m$ randomly selected entries of $M^{\sharp}$, the set of indices of which is denoted by $\Omega$, the matrix completion problem entails finding a matrix $\hat{M}$ with smallest rank that agrees with the $m$ obtained samples. That is, defining $M_{\Omega}$ to be the projection of $M$ onto the set of matrices supported on $\Omega$,
$$\hat{M}:= \text{argmin  } \text{rank}(M) \text{  s.t.  } M^{\sharp}_{\Omega}=M_{\Omega}.$$

This problem has been studied extensively in the literature \cite{candes2009exact, keshavan2009matrix,cai2010singular,candes2010power}. In general, rank-minimization is NP-hard. Therefore, one of the most common approaches is minimizing convex surrogates of the rank, such as nuclear-norm \cite{candes2010power,cai2010singular,candes2010matrix} and max-norm \cite{foygel2011concentration,cai2016matrix}. Both these convex regularizers can be used to recover the underlying matrix with as few as $O(r N\log(N))$ measurements which is close to $O(Nr)$, the number of free variables in a rank-$r$ matrix in $\mathbb{R}^{N \times N}$.\newline

A practical example where this problem arises naturally is the Netflix problem: Users and movies are arranged in the rows and columns of a matrix $M$, respectively. Each entry $M(i,j)$ represents the rating that user $i$ gives to movie $j$ \cite{koren2009matrix}. In practice it is not possible to have all the ratings from all the users, and this is where matrix completion becomes important, assuming that the true matrix is approximately low-rank. Furthermore, in most \textit{recommender systems} like Netflix, the ratings are highly quantized, sometimes even to a single bit. For example, this happens when the users just answer if they liked a movie or not. This is where the problem of 1-bit matrix completion comes into play. Here the problem is recovering a low-rank matrix from binary (1-bit) measurements of a subset of its entries.\newline

In many applications data can be represented as multi-dimensional arrays (tensors) \cite{mocks1988topographic,shashua2001linear,acar2005modeling,hu2015new}. Seismic images are ususally 3 or 5 dimensional, 2 or 4 dimensions for sources and receivers locations and 1 dimenion for time or frequency \cite{kreimer2013tensor,da2015optimization}. In hyperspectral imaging, datacubes are indexed by two spatial variables and a frequency/wavelength variable \cite{li2010tensor}. Also in some applications it is useful to introduce new information as an extra dimension. For example adding the age of users as a dimension can improve the performance of recommender systems \cite{adomavicius2011context}. The multi-dimensional structure of tensors gives them a higher expressive power compared to matrices or vectors that can be obtained by rearranging tensors.\newline

In this paper, we generalize the problem of 1-bit matrix completion to the case of tensors and analyze error bounds for recovering a low-rank tensor from partial 1-bit observations. Similar to matrices, a common assumption to make the problem feasible is assuming that the original tensor is low rank, i.e., it can be written as the sum of a few rank-one tensors. Matricizing a low-rank tensor does not increase its rank but loses the (low-dimensional) structure of the dimensions that get merged to form the rows (or columns) of the matricized tensor \cite{mu2014square}. Therefore here we consider tensor completion under M-norm (and max-qnorm) constraint, which is a robust proxy of tensor-rank \cite{ghadermarzy2017near}. In \cite{ghadermarzy2017near}, generalizing the matrix max-norm to the case of tensors led us to define two functions which are closely related, i.e., max-qnorm and atomic M-norm. To be precise, direct generalization of matrix max-norm defines a quasi-norm (max-qnorm) which is non-convex (and thus not a norm). However, we can prove near-optimal sample complexity for tensor completion using max-qnorm constrained optimization, which also gives promising results in practice. Moreover, in \cite{ghadermarzy2017near}, studying the dual of the dual of this quasi-norm resulted in defining a convex norm, atomic M-norm, which can be used to prove similar near-optimal sample complexity. In this paper, we study constrained ML estimations using both these functions. The main contributions of the paper are summarized in the next subsection.
%\subsection{Matrix completion}
%\subsection{Tensor completion}
%\subsection{1-Bit matrix completion}
\subsection{Contributions}
To the best of our knowledge, this is the first paper that analyzes the sample complexity of 1-bit tensor completion without matricization. It is worth mentioning that the closest line of research to 1-bit tensor completion is the problem of \textit{context-aware recommender systems} \cite{adomavicius2011context} and the work of \cite{karatzoglou2010multiverse}, which studies \textit{collaborative filtering} based on tensor factorization. However, in \cite{karatzoglou2010multiverse} the measurement model and the choice of regularization is different, and the sample complexity has not been analyzed.\newline

Consider a rank-$r$, order-$d$ tensor $T^{\sharp} \in \mathbb{R}^{N_1 \times \cdots \times N_d}$ where $N_i = O(N)$ for $1\leq i \leq d$. The main contributions of the current paper is as follows.
\begin{enumerate}%[label=(\roman*)]
\item We formulate and analyze the problem of 1-bit tensor completion. We show the advantage of working with a tensor directly without matricizing when $r \ll N$, both theoretically and numerically.\medskip
\item We analyze 1-bit tensor completion using M-norm constraints on the underlying tensor (this is a convex constraint). We prove that, with high probability, the mean squared error (MSE) of recovering a rank-$r$ tensor $T^{\sharp}$ from $m$ 1-bit measurements by solving an M-norm constrained log-likelihood is $O({\sqrt{r^{3d-3}}}\sqrt{\frac{Nd}{m}})$. Moreover, we analyze a related non-convex function, max-qnorm, and prove that MSE of optimizing a log-likelihood function constrained by max-qnorm is $O(\sqrt{r^{d^2-d}}\sqrt{\frac{Nd}{m}})$.\medskip
\item The M-norm gives a convex proxy for the rank of a tensor (See Remark \ref{remark_Mnorm} in Section \ref{section_Mnorm}). We derive an information-theoretic lower bound that proves the MSE of any arbitrary algorithm is $\Omega(R\sqrt{\frac{N}{m}})$ for a tensor with M-norm less than $R$. This proves that our upper bound is optimal in $N$ and the M-norm bound $R$ (but not necessarily in $r$).\medskip
\item We propose a numerical method to approximate the solution of max-qnorm constrained 1-bit tensor completion and show its advantage over 1-bit matrix completion using synthetic and real-world data. Specifically, we illustrate that one  gets significant improvement by applying the 1-bit tensor completion to context-aware recommender systems.
\end{enumerate}
%\begin{rem}
%In an accompanying paper, we analyze 1-bit tensor completion using nuclear-norm and exact-rank constraints.
%\end{rem}

\subsection{Organization}
In Section \ref{section_notations} we review some basic definitions related to tensors. Section \ref{section_prior} summarizes some related results in the literature. Section \ref{section_main_results} contains the main results regarding 1-bit tensor completion, regularized by max-qnorm and M-norm. In Section \ref{section_information} we prove a lower bound on the best error bound achievable with any arbitrary algorithm. We explain an algorithm for approximating the solution of max-qnorm constrained tensor completion in \ref{section_numerical_experiments} and present numerical results on synthetic data. We present the results on some real data in Section \ref{section_applications}. Finally in Section \ref{section_proofs} we provide all the proofs of the theorems.
\section{Notations and basics on tensors}\label{section_notations}
In this section, we introduce some basic definitions related to tensors. In particular, we introduce M-norm and max-qnorm \cite{ghadermarzy2017near} and briefly discuss some of their main properties.
\subsection{Notations}
We adopt the notations used in Kolda and Bader's review on tensor decompositions \cite{kolda2009tensor}. In what follows all universal constants are denoted by $c$ or $C$. Lowercase Greek letters are used to denote scalars, e.g., $\alpha$, $\gamma$, and $\sigma$. Matrices are denoted by $M$ or $A$ and tensors are denoted by $T$ or $X$. Other uppercase letters can either denote a matrix or a tensor depending on the context. A $d$-dimensional tensor is an element of $\bigotimes_{i=1}^{d}  \mathbb{R}^{N_i}$. The order of a tensor is the number of its dimensions (modes) and is usually denoted by $d$, and its length (size) is denoted by $(N_1, N_2,\cdots,N_d)$. When all the sizes are the same along all the dimensions, we denote $\bigotimes_{i=1}^{d}  \mathbb{R}^{N}$ as $\mathbb{R}^{N^d}$. For a tensor $X \in \bigotimes_{i=1}^{d}  \mathbb{R}^{N_i}$, elements of the tensor are alternately specified as either $X_{i_1, i_2, \cdots, i_d}$ or $ X(i_1, i_2, \cdots, i_d)$, where $1 \leq i_j \leq N_j$ for $1 \leq j \leq d$. We also use $X_{\omega}$ with $\omega=(i_1, i_2, \cdots, i_d)$, to refer to $X(i_1, i_2, \cdots, i_d)$. Inner products are denoted by $\langle \cdot , \cdot \rangle$. The symbol $\circ$ represents both matrix and vector outer products where $T=U_1 \circ U_2 \circ \cdots \circ U_d$ means $T(i_1,i_2,\cdots,i_d)=\sum_{k} U_1(i_1,k)U_2(i_2,k)$ $\cdots U_d(i_d,k)$, where $k$ ranges over the columns of the factors. In the special case when $u_j$'s are vectors, $T=u_1 \circ u_2 \circ \cdots \circ u_d$ means $T(i_1,i_2,\cdots,i_d)=u_1(i_1)u_2(i_2)\cdots u_d(i_d)$. Finally $[N]:=\{1, \cdots, N\}$.\newline

\subsection{Rank of a tensor}
A unit tensor  $U \in \bigotimes_{j=1}^{d}  \mathbb{R}^{N_j}$  is a tensor that can be written as an outer product of $d$ vectors
\begin{equation}\label{rank_one}
U=u^{(1)} \circ u^{(2)} \circ \cdots \circ u^{(d)}
\end{equation}
where $u^{(j)} \in \mathbb{R}^{N_j}$ is a unit-norm vector. The vectors $u^{(j)}$ are called the components of $U$. A rank-$1$ tensor is a non-zero scalar multiple of a unit tensor. The rank of a tensor $T$, denoted by rank($T$), is defined as the smallest number of rank-$1$ tensors that sums to $T$, i.e.,
\begin{equation*}
\text{rank}(T) = \min_r \{r: T=\sum_{i=1}^r \lambda_i U_i = \sum_{i=1}^r \lambda_i u_i^{(1)} \circ u_i^{(2)} \circ \cdots \circ u_i^{(d)}\},
\end{equation*}
where $U_i \in \onet_d$ is a unit tensor. This low-rank decomposition is also known as  CANDECOMP/PARAFAC (CP) decomposition. For a tensor $T= \sum_{i=1}^r v_i^{(1)} \circ v_i^{(2)} \circ \cdots \circ v_i^{(d)}$ we define $V_j:=[v_1^{(j)} v_2^{(j)} \cdots v_r^{(j)}]$ to be the $j$-th \textit{factor matrix} of $T$. Factor matrices can be interpreted as the higher-order generalization of collection of left (or right) singular vectors of matrices. The Frobenius norm of a tensor is defined as
\begin{equation}
\|T\|_F^2 := \langle T,T \rangle = \sum_{i_1=1}^{N_1} \sum_{i_2=1}^{N_2} \cdots \sum_{i_d=1}^{N_d} T_{i_1, i_2, \cdots, i_d}^2.
\end{equation}
\subsection{M-norm and max-qnorm}\label{section_Mnorm}
The max-qnorm of a tensor is defined as \cite{ghadermarzy2017near}
\begin{equation}\label{max-qnorm}
\|T\|_{\text{max}}:= \underset{T=U^{(1)} \circ U^{(2)} \circ \cdots \circ U^{(d)}}{\text{min}}\lbrace \prod_{j=1}^{d} \|U^{(j)}\|_{2,\infty}\rbrace,
\end{equation}
where, $\|U\|_{2,\infty} = \underset{\|x\|_2=1}{\text{sup}} \|Ux\|_{\infty}$ is the maximum row norm of $U$. In \cite{ghadermarzy2017near}, we thoroughly discussed this generalization and showed that it is a quasi-norm when $d>2$ (thus the name max-qnorm). Moreover, analyzing the dual of the dual of the max-qnorm led us to define a closely related atomic M-norm induced by the set of rank-$1$ sign tensors $T_{\pm}:=\lbrace T \in \{\pm 1\}^{N_1\times N_2 \times \cdots \times N_d}|\text{rank}(T)=1\rbrace$. Specifically, the atomic M-norm is defined as
\begin{equation}\label{M_norm}
\|T\|_{M} := \text{inf}\{t>0:T \in t\ \text{conv}(T_{\pm})\}.
\end{equation}
\noindent
\begin{rem}\label{remark_Mnorm}
The atomic M-norm is indeed a norm and is closely related to the max-qnorm \cite[Lemma 5]{ghadermarzy2017near}. The M-norm and max-qnorm are robust proxies for the rank of a tensor and both M-norm and max-qnorm of a bounded low-rank tensor is upper-bounded by a quantity that just depends on its rank and its infinity norm and is independent of $N$. In particular, assume $T \in \mathbb{R}^{N_1 \times \cdots N_d}$ is a rank-$r$ tensor with $\|T\|_{\infty} = \alpha$. Then \cite[Theorem 7]{ghadermarzy2017near}
\begin{itemize}
\item $\alpha \leq \|T\|_{M} \leq (r\sqrt{r})^{d-1} \alpha.$\medskip
\item $\alpha \leq \|T\|_{\m} \leq \sqrt{r^{d^2-d}} \alpha.$
\end{itemize}
\end{rem}
\section{Prior Work}\label{section_prior}
\subsection{1-bit matrix completion}
In the 1-bit matrix completion problem, the measurements are in the form of $Y(\omega)=\text{sign}(M^{\sharp}(\omega)+Z(\omega))$ for $\omega \in \Omega$, where $Z$ is a noise matrix. Notice that without the noise, there will be no difference in the 1-bit measurements of $M^{\sharp}$ and $\alpha M^{\sharp}$ for any positive $\alpha$. However, if we assume that the additive noise comes from a log-concave distribution, the problem can be solved by minimizing a regularized negative log-likelihood function given the measurements. In other words, the noise matrix has a dithering effect which is essential for recovery of the matrix. Under these assumptions, a nuclear-norm constrained maximum likelihood (ML) optimization was used in \cite{davenport20141} to recover $M^{\sharp}$ and it was proved that it is minimax rate-optimal under the uniform sampling model. In particular, in order to recover a rank-$r$ matrix $M^{\sharp}$ with $\|M^{\sharp}\|_{\infty} \leq \alpha$ from 1-bit measurements, they considered the weaker assumption that the matrix belongs to the set $\{ A |\  \|A\|_{\ast} \leq \alpha \sqrt{rN^2}                                                                                                                                                                                                                                                                                                                                                                                                                                                                                                                                                                                                                                                                                                                                                                                                                                                                                                                                                                                                                                                                                                                                                                                                                                                                                                                                                                                                                                                                                                                                                                                                                                                                                                                                                                                                                                                                                                                                                                                                                                                                                                                                                                                                                                                                                                                                                                                                                                                                                                                      , \|A\|_{\infty}\leq \alpha\}$ and proved
$$\frac{1}{N^2}\|M^{\sharp}-M_{\text{recovered}}\|_F^2 \leq C_{\alpha} \sqrt{\frac{rN}{m}},$$
provided that $m>CN\log(N)$. It was shown in \cite{davenport20141} that these estimates are near-optimal.  \newline

A line of research followed \cite{davenport20141} that concentrated on the recovery of low rank matrices from 1-bit measurements \cite{cai2013max,bhaskar20151,ni2016optimal}. Here we emphasize \cite{cai2013max} which considered using max-norm-constrained ML estimation instead of nuclear-norm-constrained ML. The max norm, an alternative convex proxy of the rank of a matrix, bounds the spectral norm of the rows of the low rank factors of $M^{\sharp}$, say $U$ and $V$. It is defined as $\|M^{\sharp}\|_{\m}:=\min \|U\|_{2,\infty}\|V\|_{2,\infty}\ \text{s.t.}\ M^{\sharp}=UV'$ \cite{srebro2004maximum}, where $\|U\|_{2,\infty}$ is the maximum row norm of the matrix $U$. The max norm has been extensively investigated in the machine learning community after the work of \cite{srebro2005rank} and shown to be empirically superior to nuclear-norm minimization for collaborative filtering problems. A max-norm constrained ML estimation was analyzed in \cite{cai2013max} and it was shown to be near optimal. To be precise, it was proved in \cite{cai2013max} that under some weak assumptions on the sampling distribution, with probability $1-\frac{4}{N}$,
$$\frac{1}{N^2}\|M^{\sharp}-M_{\text{recovered}}\|_F^2 \leq C_{\alpha} (\alpha\sqrt{\frac{rN}{m}}+U_{\alpha}\sqrt{\frac{\log(N)}{m}}).$$
Notice that $m=O(rN)$ is sufficient for estimating the true matrix efficiently.\newline
\subsection{1-bit tensor completion}
The generalization of 1-bit matrix completion to 1-bit tensor completion is new. However, tensor factorization for context-aware recommender systems has been investigated in \cite{karatzoglou2010multiverse, hidasi2012fast,zheng2012optimal} with two distinct differences. First they lack the 1-bit machinery introduced in \cite{davenport20141} which we show its importance in Section \ref{section_applications} and second is the lack of theoretical analysis.\newline

%Another related problem is the problem of tensor completion \cite{liu2013tensor,gandy2011tensor}. Tensor completion is the problem of recovering a low-rank tensor by observing a subset of its entries. The most relevant work is \cite{ghadermarzy2017near} which uses max-qonrm and M-norm for tensor completion and proves that $m=O(Nd)$ measurements is sufficient for efficient recovery of a rank-$r$ tensor when $r=O(1)$.
\section{Main Results}\label{section_main_results}
Next, we explain our main results. After stating the observation model and our main goal in detail in Section \ref{section_problem_formulation}, we provide recovery guarantees for 1-bit tensor completion using max-qnorm and M-norm constrained ML estimation in Section \ref{section_1bit_maxnorm}. The proofs are presented in Section \ref{section_proofs}. Necessary tools such as Rademacher complexity and discrepancy, are provided in Appendices \ref{section_appendix_radamache} and \ref{section_appendix_discrepency}.

\subsection{Problem formulation and observation model}\label{section_problem_formulation}
Given a $d$-dimensional tensor $T^{\sharp} \in \bigotimes_{i=1}^{d} \mathbb{R}^N$ (for simplicity we assume $N_i=N$) and a random subset of its indices $\Omega \subset [N] \times [N] \times \cdots [N]$, suppose that we only observe $1$-bit measurements of $T^{\sharp}$ on $\Omega$ according to the rule
\begin{equation}\label{measurements}
Y_{\omega}=
\begin{cases}
+1\ \ \ \text{with probability} \ f(T^{\sharp}_{\omega}),\\
-1\ \ \ \text{with probability} \ 1-f(T^{\sharp}_{\omega})
\end{cases}
\text{for } \omega \in \Omega.
\end{equation}
Here $f:\mathbb{R} \rightarrow [0,1]$ is a monotone differentiable function which can be interpreted as the cumulative distribution of some noise function $Z_\omega$ where the observation then becomes \cite{davenport20141}
$$
Y_{\omega}=
\begin{cases}
+1\ \ \ \text{with probability} \ T_\omega+Z_\omega \geq 0,\\
-1\ \ \ \text{with probability} \ T_\omega+Z_\omega < 0
\end{cases}
\text{for } \omega \in \Omega.$$
A standard assumption is that observed indices are chosen uniformly in random, i.e., we randomly pick $|\Omega|$ indices out of all the $N^d$ indices and then take 1-bit measurements on this set. Here, $\Omega$ is chosen at random such that $\mathbb{P}(\omega \in \Omega)=\frac{m}{N^d}$. Alternatively, the measurements may be picked uniformly at random with replacement so each index is sampled independently of the other measurements. We take a generalization of the latter model, i.e., we assume a general sampling distribution $\Pi=\{\pi_\omega\}$ for $\omega \in [N]\times[N]\times \cdots\times [N]$ such that $\sum_\omega \pi_\omega=1$ and then each index is sampled using the distribution $\Pi$. Although in this method it is possible to choose a single index multiple times, it has been proven in various contexts that sampling with replacement can be as powerful theoretically as sampling without replacement\cite{cai2016matrix,cai2013max}. Other than theoretical considerations, assuming non-uniform sampling is a better model for various applications of 1-bit matrix and tensor completion, including the Netflix problem. The challenges and benefits of using non-uniform sampling are thoroughly discussed in \cite{srebro2010collaborative,cai2013max}.\newline
\noindent

We assume that $f$ is such that the following quantities are well defined:
\begin{equation}\label{L_alpha_and_beta_alpha}
\begin{aligned}
L_{\alpha}&:= \underset{|x|\leq \alpha}{\text{sup}}  \frac{|{f'}(x)|}{f(x)(1-f(x))},\\
\beta_{\alpha}&:=\underset{|x|\leq \alpha}{\text{sup}}  \frac{f(x)(1-f(x))}{({f'}(x))^2},
\end{aligned}
\end{equation}
where $\alpha$ is the upper bound of the absolute values of entries of $T$. Here, $L_{\alpha}$ controls the \emph{steepness} of $f$, $\beta_{\alpha}$ controls its \emph{flatness}. Furthermore, we assume $f$ and $f'$ are non-zero in $[-\alpha,\alpha]$ and the quantity:
$$U_\alpha:=\underset{|x|\leq \alpha}{\text{sup}}\log(\frac{1}{f(x)(1-f(x))}),$$
is well-defined. A few well-known examples are \cite{davenport20141}:
\begin{itemize}
\item {\bf{Logistic regression/Logistic noise}} defined as $f(x)=\frac{e^x}{1+e^x}$,\newline
$L_{\alpha}=1$, $\beta_{\alpha}=\frac{(1+e^\alpha)^2}{e^\alpha}$, and $U_\alpha=2\log(e^{\frac{\alpha}{2}}+e^{\frac{-\alpha}{2}})$.

\item {\bf{Probit regression/Gaussian noise}} defined as $f(x)=\Phi(\frac{x}{\sigma})$, where $\Phi$ is the cumulative distribution function of a standard Gaussian random variable\newline
$L_{\alpha}\leq \frac{4}{\sigma}(\frac{\alpha}{\sigma}+1)$, $\beta_{\alpha}\leq \pi \sigma^2 \exp(\frac{\alpha^2}{2\sigma^2})$, and $U_\alpha\leq (\frac{\alpha}{\sigma}+1)^2$.
\end{itemize}
%L_{\alpha}:= \underset{|x|\leq \alpha}{\text{sup}}  \frac{|\dot{f}(x)|}{f(x)(1-f(x))} \ \ \ \ \text{and}\ \ \ \ \beta_{\alpha}:=\underset{|x|\leq \alpha}{\text{sup}}  \frac{f(x)(1-f(x))}{(\dot{f}(x))^2}.
%\end{equation}

The problem of $1$-bit tensor completion entails recovering $T^{\sharp}$ from the measurements $Y_{\omega}$ given the function $f$. To do this, similar to \cite{davenport20141}, we minimize the negative log-likelihood function 
\begin{equation*}
\mathcal{L}_{\Omega,Y}(X)=\sum_{\omega \in \Omega}\left(\mathds{1}_{[Y_\omega=1]}log(\frac{1}{f(X_\omega)})+\mathds{1}_{[Y_\omega=-1]}log(\frac{1}{1-f(X_\omega)})\right),
\end{equation*}
given our observation and subject to certain constraints to promote low-rank solutions.

%There are a few important quantities that should be finite for the particular choice of $F$. In particular we choose $F$ such that the following quantities are well defined:
%\begin{equation}\label{L_alpha_and_beta_alpha}

\subsection{Max-qnorm and M-norm constrained 1-bit tensor completion}\label{section_1bit_maxnorm}
In this section we analyze the problem of 1-bit tensor recovery using max-qnorm constrained ML estimation as well as M-norm constrained ML estimation. These two optimization problems are closely related, and our analysis and results are similar. Therefore, we explain the corresponding results together. To be precise, we recover $T^{\sharp}$ from the 1-bit measurements $Y_{\Omega}$, acquired following the model \eqref{measurements}, via minimizing the negative log-likelihood function
\begin{equation}\label{negative_loglikelihood}
\mathcal{L}_{\Omega,Y}(X)=\sum_{\omega \in \Omega}\left(\mathds{1}_{[Y_\omega=1]}log(\frac{1}{f(X_\omega)})+\mathds{1}_{[Y_\omega=-1]}log(\frac{1}{1-f(X_\omega)})\right),
\end{equation}
subject to $X$ having low max-qnorm or M-norm and small infinity-norm. Defining
\begin{align*}
K^T_{\m}(\alpha,R):&=\lbrace T \in \mathbb{R}^{N_1 \times N_2 \times \cdots \times N_d}: \|T\|_{\infty} \leq \alpha, \|T\|_{\m} \leq R \rbrace,\\
K^T_{M}(\alpha,R):&=\lbrace T \in \mathbb{R}^{N_1 \times N_2 \times \cdots \times N_d}: \|T\|_{\infty} \leq \alpha, \|T\|_{M} \leq R \rbrace,
\end{align*}
we analyze the following optimization problems:
\begin{align}
\hat{T}_{\m} &= \underset{X}{\text{arg min }} \mathcal{L}_{\Omega,Y}(X) \ \ \ \ \text{subject to}\ \ \ \  X \in K^T_{\m}(\alpha,R_{\m}). \label{optimization_maxnorm} \\
\hat{T}_{M} &= \underset{X}{\text{arg min }} \mathcal{L}_{\Omega,Y}(X) \ \ \ \ \text{subject to}\ \ \ \  X \in K^T_{M}(\alpha,R_M). \label{optimization_Mnorm} 
\end{align}
\begin{rem}\label{remark_flat_tensors}
Flat low-rank tensors are contained in $K^T_{\m}(\alpha,R_{\m})$ and $K^T_{M}(\alpha,R_M)$ ($R_{\m}$ and $R_{M}$ are quantities that just depend on $\alpha$ and the rank of the tensor) \cite[Theorem 7]{ghadermarzy2017near}. M-norm, max-qnorm and infinity norm are all continuous functions and therefore, both $K^T_{\m}(\alpha,R_{\m})$ and $K^T_{M}(\alpha,R_M)$ contain approximately low-rank tensors. We assume the data is generated from one of these sets of approximate low-rank tensors.
\end{rem}
\begin{thm}\label{theorem_maxnorm}
Suppose that we have $m$ 1-bit measurements of a subset of the entries of a tensor $T^{\sharp} \in \mathbb{R}^{N^d}$ following the probability distribution explained above where $f$ is chosen such that $L_\alpha$, $\beta_\alpha$, and $U_\alpha$ are well defined. Moreover, the indices of the tensor are sampled according to a probability distribution $\Pi$. Then assuming $\|T^{\sharp}\|_{\m} \leq R_{\m}$, there exist absolute constants $C_{\m}$ and $C_M$ such that for a sample size $2\leq m \leq N^d$ and for any $\delta>0$ the maximizer $\hat{T}_{\m}$ of \eqref{optimization_maxnorm} satisfies:
\begin{equation}\label{upper_bound}
\|T^{\sharp} - \hat{T}_{\m}\|_{\Pi}^2:=\sum_\omega \pi_\omega (\hat{T}_{\m}(\omega)-T^{\sharp}(\omega))^2 \leq C_{\m} c_2^d \beta_\alpha \lbrace L_\alpha R_{\m} \sqrt{\frac{dN}{m}}+U_\alpha \sqrt{\frac{\log(\frac{4}{\delta})}{m}}\rbrace
\end{equation}
with probability at least $1-\delta$. Similarly, assuming $\|T^{\sharp}\|_{M} \leq R_{M}$, for any $\delta>0$, the maximizer $\hat{T}_{M}$ of \eqref{optimization_Mnorm} satisfies:
\begin{equation}\label{upper_bound_Mnorm}
\|T^{\sharp} - \hat{T}_{M}\|_{\Pi}^2=\sum_\omega \pi_\omega (\hat{T}_{M}(\omega)-T^{\sharp}(\omega))^2 \leq C_M \beta_\alpha \lbrace L_\alpha R_M \sqrt{\frac{dN}{m}}+U_\alpha \sqrt{\frac{\log(\frac{4}{\delta})}{m}}\rbrace
\end{equation}
with probability at least $1-\delta$.
\end{thm}
\begin{cor}
A rank-$r$ tensor $T^{\sharp}_r$ with $\|T^{\sharp}_r\|_{\infty} \leq \alpha$ satisfies $\|T^{\sharp}_r\|_{\m} \leq \sqrt{r^{d^2-2}} \alpha$ and therefore, the maximizer $\hat{T}_{\m}$ of \eqref{optimization_maxnorm} with $R_{\m}=\sqrt{r^{d^2-2}} \alpha$ satisfies
$$ \|T^{\sharp}_r - \hat{T}_{\m}\|_{\Pi}^2 \leq C_{\m} c_2^d \beta_\alpha \lbrace L_\alpha \sqrt{r^{d^2-2}} \alpha \sqrt{\frac{dN}{m}}+U_\alpha \sqrt{\frac{\log(\frac{4}{\delta})}{m}}\rbrace,$$
with probability greater than $1-\delta$. Moreover, $\|T^{\sharp}_r\|_{M} \leq (r\sqrt{r})^{d-1} \alpha$ and the maximizer $\hat{T}_{M}$ of \eqref{optimization_Mnorm} with $R_{M}=(r\sqrt{r})^{d-1} \alpha$ satisfies
$$\|T^{\sharp}_r - \hat{T}_{M}\|_{\Pi}^2 \leq C_M \beta_\alpha \lbrace L_\alpha (r\sqrt{r})^{d-1}  \sqrt{\frac{dN}{m}}+U_\alpha \sqrt{\frac{\log(\frac{4}{\delta})}{m}}\rbrace,$$
with probability greater than $1-\delta$.
\end{cor}
\begin{rem}
There are two differences in the right hand sides of \eqref{upper_bound} and \eqref{upper_bound_Mnorm}. First is the difference in $R_{\m}$ and $R_M$ which is due to the different upper bounds on max-qnorm and M-norm of a bounded rank-$r$ tensor. The second difference is in the constants $C_{\m} c_2^d$ and $C_M$ where $C_M < \frac{C_{\m}}{c_1}$ where $c_1$ and $c_2$ are small constant derived from the generalized Grothendieck's theorem \cite[Theorem 13]{ghadermarzy2017near},\cite{tonge1978neumann} ($c_1 <0.9$ and $c_2<\sqrt{2}$). Note for a tensor $T$, $\|T\|_{M} \leq c_1 c_2^d \|T\|_{\m}$ and $C_M < \frac{C_{\m}}{c_1}$ which implies that \eqref{upper_bound_Mnorm} is tighter than \eqref{upper_bound}. Moreover, for $d \geq 2$, the M-norm bound of a low rank tensor is smaller than its max-qnorm bound \cite{ghadermarzy2017near}.
\end{rem}
\begin{rem}
Assuming that all entries have a nonzero probability of being observed, i.e., assuming that there exist a constant $\eta \geq 1$ such that $\pi_\omega \geq \frac{1}{\eta N^d}$, for all $\omega \in [N_1] \times \cdots \times [N_d]$, we can simplify \eqref{upper_bound} as
$$\frac{1}{N^d} \|T^{\sharp} - \hat{T}_{\m}\|_{F}^2 \leq C_{\eta} \beta_\alpha \lbrace L_\alpha R_{\m} \sqrt{\frac{dN}{m}}+U_\alpha \sqrt{\frac{\log(dN)}{m}}\rbrace.$$
A similar bound can be obtained for M-norm constrained ML estimation \eqref{upper_bound_Mnorm} as well.
\end{rem}
\begin{rem}
The terms $r^{3d-3}$ and $r^{d^2-d}$ in the upper bounds come from the M-norm and max-qnorm of low rank tensors. We believe this is sub-optimal in $r$. Indeed, when $d=2$, we know that both these upper bound can be improved to $\sqrt{r}$ instead of $r^{\frac{3}{2}}$ and $r^2$ .
\end{rem} 
\begin{rem}
For a rank-$r$ tensor $T^{\sharp}$ with $\|T^{\sharp}\|_{\infty}\leq \alpha$, it is sufficient to choose $R_{\m}=r^{\frac{d^2-d}{2}} \alpha$ in \eqref{optimization_maxnorm} and $R_M = r^{\frac{3d-3}{2}} \alpha$ in \eqref{optimization_Mnorm} for $T^{\sharp}$ to be feasible \cite{ghadermarzy2017near}. This proves that efficient recovery can be achieved when $m=O(r^{d^2-d} N)$ using \eqref{optimization_maxnorm} and $m=O(r^{3d-3} N)$ using \eqref{optimization_Mnorm}. These are significant improvements over matricizing, which would require $\Omega(r N^{\frac{d}{2}})$ measurements when $r \ll N$.
\end{rem}
 
\begin{rem}
1-bit tensor completion can recover the magnitude of the tensor because of the dithering effect of the noise we add before sampling. If the SNR is too low, we might end up estimating the noise by a low rank tensor and if the SNR is too high, we risk loosing magnitude information. In Section \ref{section_numerical_experiments}, we conduct some experiments and discuss the optimal noise level using synthetic data.
\end{rem}
 
\begin{rem}
1-bit matrix completion using max-norm is analyzed in \cite{cai2013max}. A direct generalization resulted in the error bound \eqref{upper_bound} which recovers their result when $d=2$. However, the direct generalization of max-norm to tensors is a quasi-norm and the resulting constraint is non-convex. The closely related upper bound using M-norm \eqref{upper_bound_Mnorm} is governed by a convex optimization and gives similar results when $d=2$. In this case the only difference is the Grothendieck's constant $1.67<c_1 c_2^2<1.79$ \cite{jameson1987summing}.
\end{rem}

\section{Information theoretic lower bound}\label{section_information}
We use a classical information-theoretic technique that shows that with limited amount of information we can distinguish only a limited number of tensors from each other. To that end, as in \cite{davenport20141,cai2013max}, we first construct a set of tensors $\chi$ such that for any two distinct members $X^i$ and $X^j$ of $\chi$ $\|X^i - X^j \|_F$ is large. Therefore, we should be able to distinguish the tensors in this set if the recovery error of an arbitrary algorithm is small enough. However, Fano's inequality will imply that the probability of choosing the right tensor among this set is going to be small and therefore force a lower bound on the performance of any arbitrary recovery algorithm.\newline

%
%In this section we provide a lower bound for the reconstruction error when recovering tensors from 1-bit random measurements. In particular, by using a similar information theoretic argument as the one in \cite{davenport20141}, we prove that on a particular set of tensors, any arbitrary algorithm will not be able to recover a close enough approximation of all the tensors in that set. The same approach was taken in both \cite{davenport20141,cai2013max} and therefore, we will not go into the details of the proofs and refer to those papers for more detailed explanation of these type of methods. However, to be more precise, we use a classical information-theoretic technique which shows that with limited amount of information we can distinguish only a limited number of tensors from each other. To that end, we first construct a set of tensors $\chi$ such that for any two distinct members $X^i$ and $X^j$ of $\chi$ $\|X^i - X^j \|_F$ is large. Therefore, we should be able to distinguish the tensors in this set if the recovery error of an arbitrary algorithm is small enough. However, Fano's inequality will imply that the probability of choosing the right tensor among this set is going to be small and therefore force a lower bound on the performance of any arbitrary recovery algorithm.\newline

The lower bound is achieved using Lemma 18 of \cite{ghadermarzy2017near}, which constructs a packing set for the set of low M-norm tensors given by
$$K^T_{M}(\alpha,R_M):=\lbrace T \in \mathbb{R}^{N_1 \times N_2 \times \cdots \times N_d}: \|T\|_{\infty} \leq \alpha, \|T\|_{M} \leq R_M \rbrace.$$
%The constructed set can be also used as a packing set for bounded low rank tensors given by
%$$C_{\text{rank}}(r,\alpha) := \{ T \in \bigotimes_{i=1}^d \mathbb{R}^N : \|T\|_{\infty} \leq \alpha, \text{rank}(T) \leq r \}.$$
We include Lemma 18 of \cite{ghadermarzy2017near} and an adoptation of the proof in \cite{davenport20141} in Section \ref{section_proof_theorem_upperbound}. In the following $1.68 < K_G:=c_1 c_2^2 <1.79$ is the Grothendieck's constant.% In the following theorem $K$ can be either $K^T_{M}(\alpha,R)$ or $C_{\text{rank}}(r,\alpha)$.

\begin{thm}\label{theorem_information_theoritic}
Fix $\alpha$, $R_M$, and $N$ and set $r:=\floor{(\frac{R_M}{\alpha K_G})^2}$ such that $\alpha \geq 1$, $\alpha^2 r N \geq C_0$, $r \geq 4$ and $r \leq O(\frac{N}{\alpha^2})$. Let $\beta_{\alpha}$ be defined as in \eqref{L_alpha_and_beta_alpha}. Suppose $f^{'}(x)$ is decreasing for $x>0$ and for any tensor $T$ in the set $K_M^T(\alpha, R_M)$, assume we obtained $m$ measurements (denoted by $Y_{\Omega}$) from a random subset $\Omega$ of the entries  following the model \eqref{measurements}. Consider any arbitrary algorithm that takes these measurements and returns an approximation $\hat{T}(Y_{\Omega})$. Then there exist a tensor $T \in K^T_{M}(\alpha,R_M)$ such that
$$\mathbb{P} \left(\frac{1}{N^d} \|T-\hat{T}(Y_{\Omega})\|_F^2 \geq \text{min} \big(C_1, C_2 \sqrt{\beta_{\frac{3\alpha}{4}} }\frac{R_M}{K_G}\sqrt{\frac{N}{m}} \big)\right) \geq \frac{3}{4}$$
%Similarly, there exist a tensor $T \in C_{\text{rank}}(r,\alpha) $ such that
%$$\mathbb{P} \left(\frac{1}{N^d} \|T-\hat{T}(Y_{\Omega})\|_F^2 \geq \text{min} \big(C_1, C_2 \alpha \sqrt{\beta_{\frac{3\alpha}{4}} }\sqrt{\frac{rN}{m}} \big)\right) \geq \frac{3}{4}.$$
\end{thm}
\begin{rem}
Theorems \ref{theorem_maxnorm} and \ref{theorem_information_theoritic} prove that the sample complexity achieved by \eqref{optimization_Mnorm} is optimal in both $N$ and $R_M$. Evidence in the matrix case suggests that the recovery lower bound for the set of low rank tensors (instead of low M-norm) should be bounded by $\frac{N}{m}$ instead of $\sqrt{\frac{N}{m}}$ provided that $m$ is large enough. We postpone a lower bound for exact low-rank tensors to future work.
\end{rem}

\section{Numerical results}\label{section_numerical_experiments}
In this section, we present some numerical results of using \eqref{optimization_maxnorm} for 1-bit tensor completion. The optimization problem we want to solve is in the form of
\begin{equation}\label{optimization_general}
\underset{X}{\text{arg min }} \mathcal{L}_{\Omega,Y}(X) \ \ \ \ \text{subject to}\ \ \ \  X \in K^T_{\m}(\alpha,R_{\m}),
\end{equation}
where $\mathcal{L}_{\Omega,Y}(X)$ is the negative log-likelihood function defined in \eqref{negative_loglikelihood}. Similar to the matrix-completion problem $\mathcal{L}_{\Omega,Y}$ is a smooth function from $\mathbb{R}^{N^d} \rightarrow \mathbb{R}$. First, we mention a few practical considerations.
\subsection*{Max-qnorm constrained 1-bit tensor completion}
We are not aware of any algorithm that is able to (approximately) solve the M-norm constrained problem \eqref{optimization_Mnorm}. Therefore, instead, we focus on algorithms to (approximately) solve the max-qnorm constrained 1-bit tensor completion problem. To this end, we employ a similar approach to the one in \cite{ghadermarzy2017near} and use the low rank factors to estimate the max-qnorm of a tensor. In particular, defining $f(V_1, \cdots, V_d):=\mathcal{L}_{\Omega,Y}(V_1 \circ \cdots \circ V_d)$ the problem becomes
\begin{equation}\label{optimization_maxnorm_algorithm}
\text{min}\ f(V_1, \cdots, V_d) \text{ subject to } \max_i (\|V_i\|_{2,\infty}) \leq R_{\m}^{\frac{1}{d}}\ , \|V_1 \circ \cdots \circ V_d\|_{\infty} \leq \alpha.
\end{equation}
As explained in details in \cite{ghadermarzy2017near}, in the definition of max-qnorm (\ref{max-qnorm}), there is no limit on the size of the low rank factors. Due to computational restrictions we limit the factor size by twice the dimension size, i.e., $V_{i} \in \mathbb{R}^{N \times 2N}$. Although we end up approximating the solution of \eqref{optimization_maxnorm_algorithm}, our numerical results shows that, for synthetic tensors with rank smaller than $N$, $2N$ is a large enough factor size since results do not improve significantly when larger factors are used, at least in our experiments where $N<100$. Another practical benefit of \eqref{optimization_maxnorm_algorithm} is that we can use alternating minimization on the small factors instead of the whole tensor, which is crucial in applications where storing the whole tensor is too expensive. We drop the infinity-norm constraint in our simulations-See Remark \ref{remark_infconstraint} below.\newline

The next practical problem is in choosing $R_{\m}$. Although theory suggests an upperbound of $r^{\frac{d^2-d}{2}}\alpha$ \cite{ghadermarzy2017near}, in practice this bound can be much higher than the actual max-qnorm of a rank-$r$ tensor. Moreover, in many practical applications, we do not have an accurate upper bound on the rank of the tensor (or the rank of a good approximation of a tensor). Therefore, we use cross validation to estimate the optimal value of $R_{\m}$.\newline

In summary, to approximate the solution of \eqref{optimization_maxnorm_algorithm}, in our experiments we use 10\% of the available data for cross validating the choice of $R_{\m}$ and optimize over the factors $V_1$ to $V_d$ alternatingly while at each iteration, defining $f_i(V_i) := f(V_1, \cdots,V_i, \cdots, V_d)$, we solve the simpler sub-problem
\begin{equation}\label{optimization_subproblem}
\underset{V_i}{\text{min}}\ f_i(V_i) \text{ subject to }  \|V_i\|_{2,\infty} \leq R_{\m}^{\frac{1}{d}},
\end{equation}
while fixing all the other factors $V_j$, $j \neq i$. We solve each sub-problem \eqref{optimization_subproblem} by employing a projected gradient algorithm. The algorithm updates all the factors
\begin{equation}\label{equation_algorithm}
[V_i] \leftarrow \mathbb{P}_{C}([V_i - \gamma \bigtriangledown f_i]).
\end{equation}
where, $\mathbb{P}_C$ simply projects the factor onto the set of matrices with $\ell_{2,\infty}$ norm less than $R^{\frac{1}{d}}$. This projection looks at each row of the matrix and if the norm of a row is bigger than $R^{\frac{1}{d}}$, it scales that row back down to $R^{\frac{1}{d}}$ and leaves other rows unchanged.\newline
\begin{rem}\label{remark_infconstraint}
The constraint on the boundedness of the tensor ($\|V_1 \circ \cdots \circ V_d\|_{\infty} \leq \alpha$) can also be incorporated into \eqref{optimization_subproblem}, which introduces new challenges to the algorithm. First, the exact projection onto the set of $\{V_i | \|V_1 \circ \cdots \circ V_d\|_{\infty} \leq \alpha\}$ is not as straightforward as the projection onto $\{V_i | \|V_i\|_{2,\infty} \leq R_{\m}^{\frac{1}{d}}\}$. An approximate projection by rescaling the factor via
$$V_i = \frac{\alpha}{\|V_1 \circ \cdots \circ V_d\|_{\infty}} V_i,\ \text{if } \|V_1 \circ \cdots \circ V_d\|_{\infty}>\alpha $$
 was introduced in \cite{cai2013max}. On the other hand, the exact projection can be formulated as a quadratic linear program, which is very computationally expensive. The second complication is that adding this constraint makes the constraint set $C$ the intersection of two convex sets, which again makes exact projection expensive.\newline
In our synthetic experiments, we noticed that adding this constraint does  not change the results significantly, especially for low-rank tensors. Furthermore, in our applications, we concentrate on the performance of the algorithm in recovering the sign of the tensor which reduces the importance of projecting onto the set of bounded tensors and therefore, in this section we just report the results obtained by (approximately) solving \eqref{optimization_subproblem}.
\end{rem}

\begin{rem}
The optimization problem \eqref{optimization_subproblem} can be solved in parallel over the rows of $V_i$ as both the objective function and the regularizer ($\| . \|_{2,\infty}$) are decomposable in the rows.
\end{rem}

\subsection{Synthetic experiments}\label{section_synthetic}
\begin{figure}[t]
\centering
\includegraphics[width=1.1\textwidth]{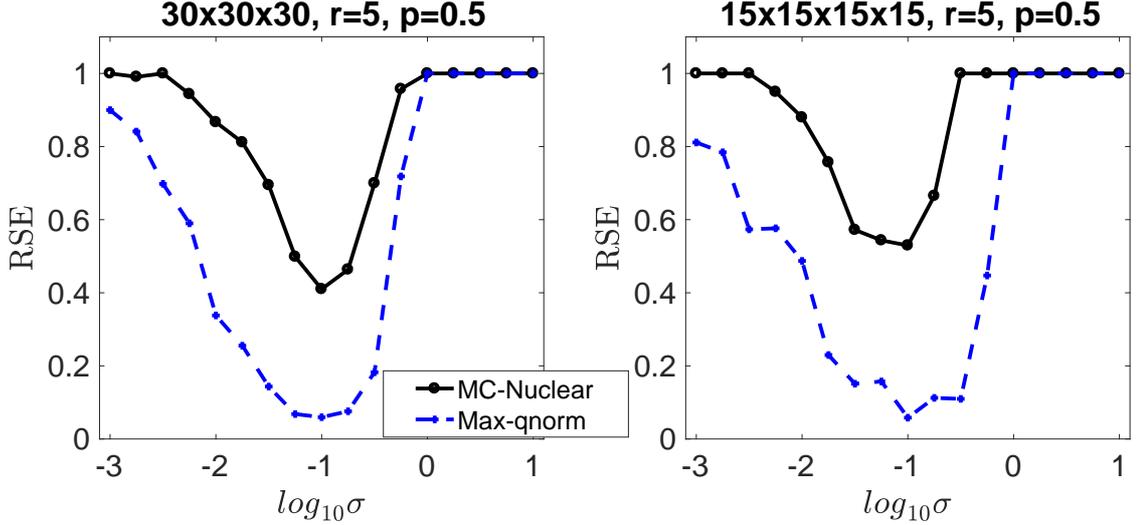}
\caption{Average relative error obtained from nuclear-norm constrained matricization and max-qnorm. The results show the relative error over a range of different value of $\sigma$ averaged for 10 random rank-$5$ tensors. From left to right the original tensor is 2 and 3 dimensional. Results show that the noise variance should not be too small or too big.}
\label{1bitTCsigmar5}
\end{figure}
In this section we present extensive numerical results on synthetic experiments to compare the performance of 1-bit tensor completion using max-qnorm with matricization using nuclear norm \cite{davenport20141}. An important component of the measurements which effects the results is the choice of the function $f$. We investigate the optimal choice of $\sigma$ when we use the Gaussian noise, i.e., $f(x)=\Phi(\frac{x}{\sigma})$ in the log-likelihood function \eqref{negative_loglikelihood}. Figure \ref{1bitTCsigmar5} shows the results for $d=3,4$, where $n=30,15$ respectively. Here, we average the recovery error over 20 realizations when $r=5$ and $\sigma$ varies from $0.001$ to $10$. Moreover, in all cases, $m=|\Omega| = \frac{N^d}{2}$.\newline

In Figure \ref{1bitTCsigmar5}, we show the results for matricized nuclear-norm constrained 1-bit matrix completion (MC-Nuclear) and max-qnorm as explained above. The underlying low-rank tensor is generated from low rank $N \times r$ factors with entries drawn i.i.d. from a uniform distribution on $[-1,1]$. The tensor is then scaled to have $\|T^{\sharp}\|_{\infty}=1$. The figure shows the relative squared error 
$$RSE:=\frac{\|T_{\text{recovered}}-T^{\sharp}\|_F^2}{\|T^{\sharp}\|_F^2}$$
against different values of $m$. The MC-Nuclear results show the case where we matricize the tensor first and then solve a 1-bit matrix completion problem using the algorithm in \cite{davenport20141}, but with added cross validation for the optimal nuclear-norm bound.  Similar to the matrix case \cite{davenport20141}, the results deteriorate in the high-noise regime as we end up measuring the noise, which in the 1-bit regime looks like a coin toss. Choosing small noise on the other limits the dithering effect which is essential for successful recovery. Notice that the same experiment can be done with various values of $r$, and $m$. However, finding the optimal value of $\sigma$ for different cases is time consuming and unnecessary as for most cases $\sigma=0.1$ seems to be close to the optimal value when $\|T^{\sharp}\|_{\infty} =1$.\newline

\begin{figure}[t]
\centering
\includegraphics[width=1.1\textwidth]{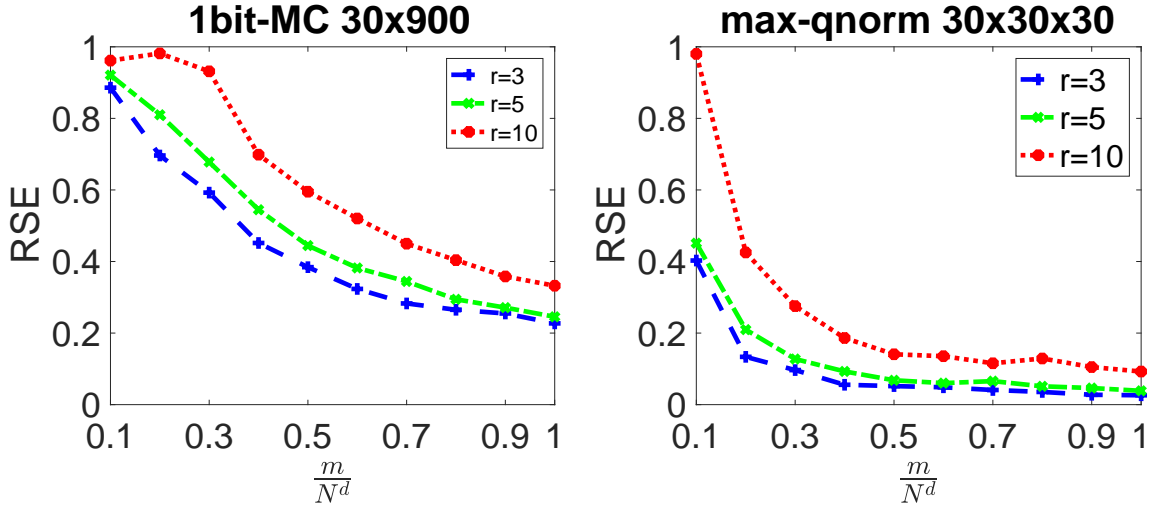}
\caption{Average relative error for recovering a $30 \times 30 \times 30$ tensor with various ranks, using $m$ 1-bit measurements sample for a range of different value of $m$. The left figure shows the results of matricizing and the right figure shows the results using \eqref{equation_algorithm}.}
\label{1bitTC_d3_N30}
\end{figure}

Once we fix the choice of $\sigma$ we move on to comparing results for different ranks and tensor sizes. We present results for various rank $r \in \{3,5,10\}$ and various sample sizes $ \frac{m}{N^d} \in \{ 0.1, 0.2, \cdots, 1\}$. Figure \ref{1bitTC_d3_N30} shows the results for 3-dimensional tensors $T^{\sharp} \in \mathbb{R}^{30 \times 30 \times 30}$. A balanced matricization is not possible in this case and therefore, tensor completion has an added advantage when the tensor has an odd order. We can see this advantage clearly in Figure \ref{1bitTC_d3_N30}. When the number of observations is small, e.g., when $\frac{m}{N^d}=0.3$ the average error of 1-bit tensor completion is four to six times better than matricizing. Note that in both matricization and max-qnorm constrained recovery we use 10\% of the observations to validate the choice of the regularization parameter.\newline

\begin{figure}[h]
\centering
\includegraphics[width=1.1\textwidth]{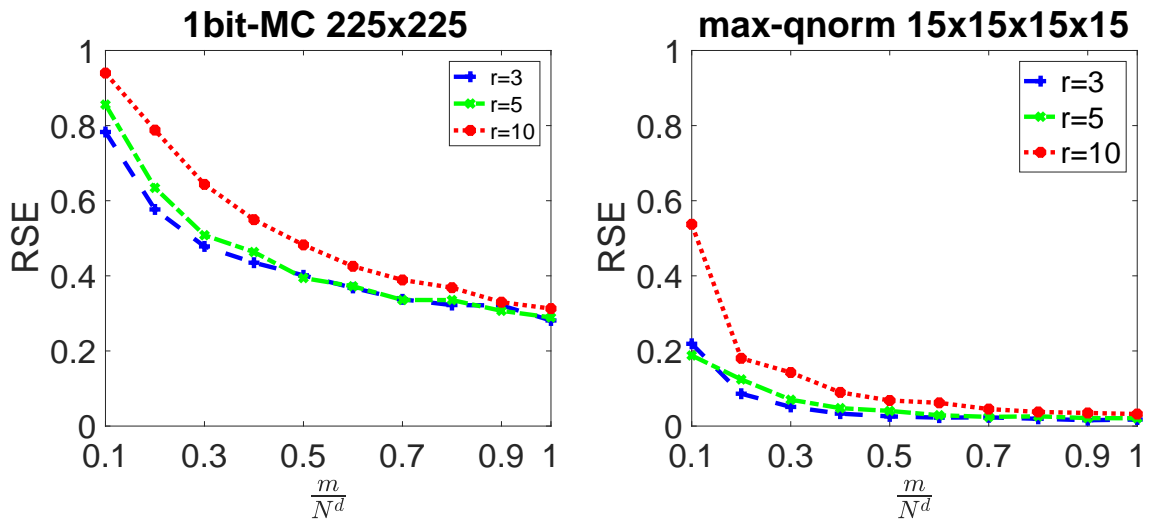}
\caption{Average relative error for recovering a $15 \times 15 \times 15 \times 15$ tensor with various ranks, using $m$ 1-bit measurements sample for a range of different value of $m$.}
\label{1bitTC_d4_N15}
\end{figure}

Figure \ref{1bitTC_d4_N15} shows the results for $T^{\sharp} \in \mathbb{R}^{15 \times 15 \times 15 \times 15}$. Again, we see a significant improvement over matricization.

\section{Applications}\label{section_applications}
In this section we test the performance of \eqref{optimization_maxnorm_algorithm} on real-world data, in particular, for predicting exam scores and user ratings. We show the improvements gained by solving 1-bit tensor completion instead of matricizing. We transform the ratings (or scores) to 1-bit data by replacing them with whether they are above or below some approximate mean-rating. Next we investigate the ability of 1-bit tensor completion for predicting if an observed rating is above or below average. In all the applications we know the maximum value of the ratings and use this value to rescale the data to have unit infinity norm. This helps us in choosing the the appropriate function $f$ based on our synthetic experiments in Section \ref{section_synthetic}, i.e., Gaussian noise with $\sigma=0.1$. In Remark \ref{remark_sigma_value} we explained the consequences of this choice in more details.\newline

To be more precise, we explain the general application setup we use in this section briefly. Assume $T^{\sharp}$ is the true tensor and we have access to a subset of its entries for training ($T^{\sharp}_{\Omega_{\text{train}}}$) and another subset of the entries for testing ($T^{\sharp}_{\Omega_{\text{test}}}$). $\Omega_{\text{train}}$ and $\Omega_{\text{test}}$ are two disjoint subsets of the entries. Below is a brief recipe of the experiments done in this section.
\begin{enumerate}
\item Scale the observed entries to have unit infinity norm.\medskip
\item Using an approximate mean-value $\eta$, take $\text{sign}(T^{\sharp}_{\Omega_{\text{train}}} - \eta)$ to be the 1-bit measurements. Notice that with synthetic low rank tensors we add a dithering noise using some function $f$ before taking the sign (see \eqref{measurements}). However, our experiments showed that assuming the noise to be intrinsic in the data gives slightly better results than adding. The same approach was taken in \cite{davenport20141}.\medskip
\item Use \eqref{optimization_maxnorm_algorithm} to get an initial estimate $\hat{T}_{\text{init}}$.\medskip
\item Add the approximate mean-value $\eta$ and scale the resulting tensor back to get the final $\hat{T}$, which is an approximation of $T^{\sharp}$.\medskip
\item Evaluate $\hat{T}$ by comparing $\hat{T}$ and $T^{\sharp}$ on the test indices.
\end{enumerate}

\begin{rem}\label{remark_sigma_value}
In the applications, we empirically observe that we get better predictions if we do not dither the original tensor, i.e., if we assume the noise to be implicit in the data. Moreover, as explained in Remark \ref{remark_infconstraint}, in the applications section we also ignore the infinity-norm constraint. Therefore, in theory changing the value of $\sigma$ should not change the results of predicting the sign of the tensor. For example, if we consider recovering the original tensor with two noise functions $f_1(x)=\Phi(\frac{x}{\sigma_1})$ and $f_2(x)=\Phi(\frac{x}{\sigma_2})$, we would have
$$\Phi(\frac{T(\omega)}{\sigma_1}) = \Phi(\frac{\frac{\sigma_2}{\sigma_1} T(\omega)}{\sigma_2}),$$
which shows that the recovered tensors achieved by $f_1$ and $f_2$ should be scalar multiples of each other and the sign predictions should be the same. In Figure \ref{figure_scalability} we illustrate this by showing the results of applying noise functions with different values of $\sigma$ for recovering a rank-5 tensor in $\mathbb{R}^{30 \times 30 \times 30}$ whose 1-bit measurement have been obtained with $\sigma=0.15$. The plot in the left shows the average RSE and as expected the best recovery of the original tensor is obtained when we use the true noise function that we originally used to get the observations. However, the right plot shows that the percentage of correct sign predictions is very robust to the choice of $\sigma$ and all the results are very close to each other which supports the above discussion.
\begin{figure}[t]
\centering
\includegraphics[width=1.1\textwidth]{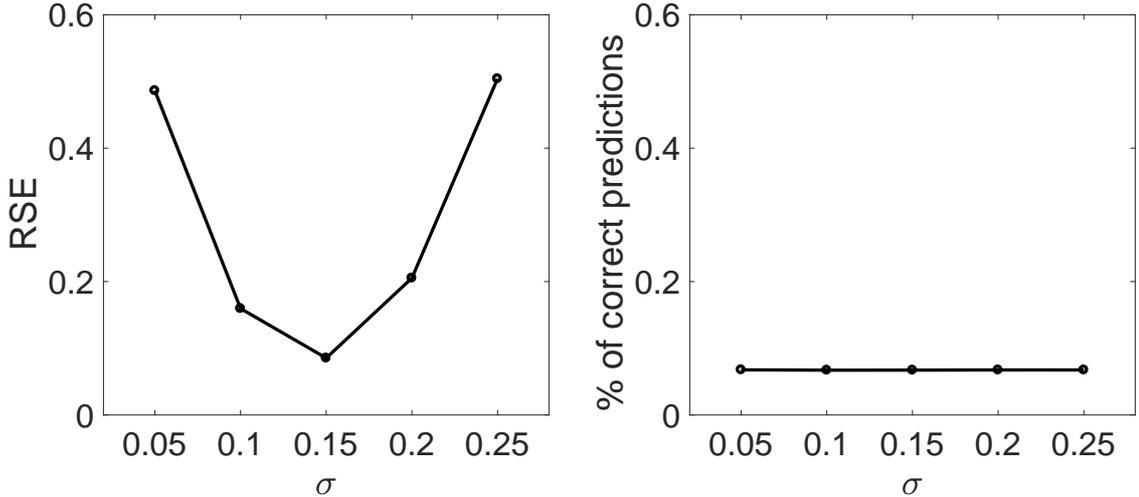}
\caption{The left figure shows the average relative error of recovering a rank-$5$, $30 \times 30 \times 30$ tensor with different values of $\sigma$ when the 1-bit measurements were obtained by using $\sigma=0.15$. The right figure shows the percentage of correct sign predictions.}
\label{figure_scalability}
\end{figure}

\end{rem}
\subsection{Score prediction}
\begin{figure}
\includegraphics[scale=0.35]{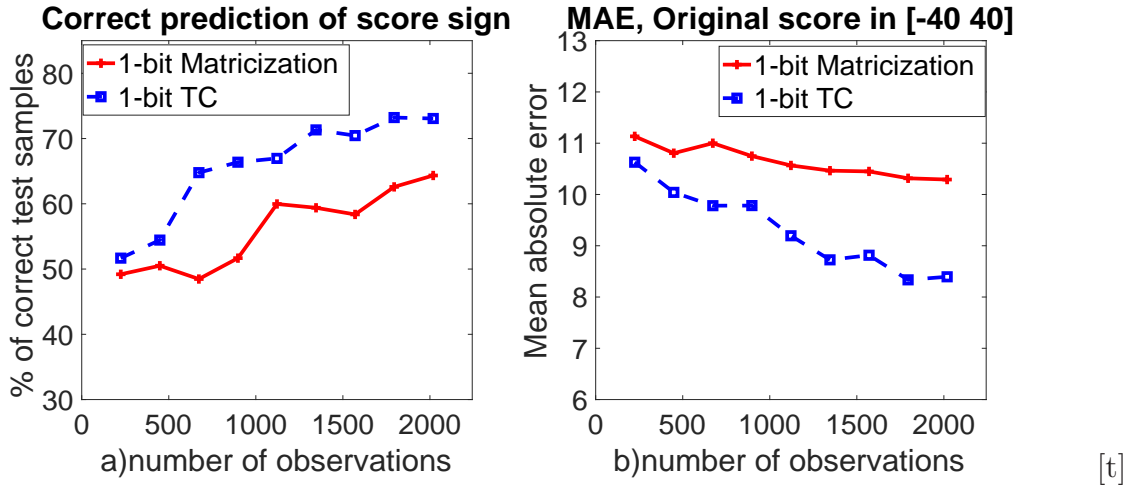}[t]
\caption{Results of applying 1-bit matrix and tensor completion to partially observed verbal scores tensor \cite{mcgrath1986british} to determine whether unobserved scores are above or below average. The left figure shows the percentage of correct predictions and the right figure shows the absolute mean error. The scores are in the range [-40 40]}
\label{2lev_scores}
\end{figure}

In this section, we apply our algorithm to the data on pupil attainments in schools in Scotland which contains the information of 3,435 children who attended 148 primary schools and 19 secondary schools in Scotland\footnote{Available at http://www.bristol.ac.uk/cmm/learning/mmsoftware/data-rev.html}. We generate a 5-dimensional tensor in $\mathbb{R}^{148 \times 19 \times 2 \times 4 \times 10}$ which includes the information of primary school, secondary school, gender, social class and attainment score of the students and estimate the verbal reasoning score of the students based on varying number of 1-bit information (from 230 samples to 2100 samples which is less than 1$\%$ of the total entries of the tensor). The scores are ranged from -40 to 40 and we take the 1-bit information to be the sign of the score. We use 10\% of the observations for cross validating the choice of parameters and another 10\% of the total scores as a test set. Figure \ref{2lev_scores}.a shows the percentage of correct sign predictions and Figure \ref{2lev_scores}.b shows the mean absolute error
$$MAE:=\frac{\sum_{\omega \in \Omega_t}|T^{\sharp}(\omega)-\hat{T}(\omega)|}{|\Omega_t|}$$
on the test set. Notice that the scores are in the range of -40 to 40 and the Mean absolute error is in the range of 8 to 10. The matrix completion results refer to matricizing the tensor to a $592 \times 380$ matrix by putting the first and fourth dimension in the rows and the rest in the columns. This matricization results in the most balanced rearrangement which is the recommended way for better results \cite{mu2014square}. In both figures we can see that using tensor completion outperforms matrix completion significantly.

\subsection{In-car music recommender system}
In this section we apply our algorithm to an in-car music recommender system \cite{baltrunas2011incarmusic} which contains 4000 ratings from 42 users and 140 music tracks. The ratings are acquired via a mobile application and each rating is accompanied with one of the 26 context information which ranges from the landscape to the weather condition or the driver's mood. The resulting tensor is of size $42 \times 140 \times 26$ and we use 3200 ratings for the minimization (less than 2$\%$ of the entries of the tensor), 400 ratings for validating the choice of max-qnorm bound and the other 400 ratings as test data. Table \ref{table:incarmusic} shows the results using tensor completion while considering the context information and matrix completion. Both matrix and tensor completion results are obtained by the max-qnorm constrained tensor completion algorithm explained in \cite{ghadermarzy2017near}. The table shows the correct 1-bit predictions made by each algorithm considering their original ratings in the first six columns and the total average in the last row. Bringing context into the reconstruction results in an impressive improvement of at least 17\% and moreover, using 1-bit information does not change the results too much compared to using the full information. Note that the results are averaged over 10 different random training and test sets.
\begin{table}[h]
\centering
\begin{tabular}{|r||c | c | c | c | c | c || c|}
\hline
Original rating & 0 & 1 & 2 & 3 & 4 & 5 & Overall \\ \hline \hline
1-bit matrix completion & 75\% & 65\% & 54\% & 51\% & 70\% & 65\% & 60\% \\ \hline
multi-bit matrix completion & 76\% & 62\% & 50\% & 54\% & 53\% & 61\% & 57\% \\ \hline
1-bit TC (context-aware) & \bf{80}\% & \bf{89}\% & 58\% & \bf{65}\% & \bf{78}\% & 85\% & \bf{77}\% \\ \hline
multi-bit TC(context-aware) & \bf{80}\% & 86\% & \bf{60}\% & 64\% & 77\% & \bf{90}\% & 76\% \\ \hline
\end{tabular}
\caption{\small \sl Results of a comparison between 1-bit and multi-bit matrix and tensor completion algorithms on incar-music data \cite{baltrunas2011incarmusic} for predicting whether the unobserved ratings were above or below average. The multi-bit methods uses the original ratings from 0 to 5 and the context-aware methods include the context information such as time of the day, weather, location of driving and mood.}
 \label{table:incarmusic}
\end{table}
\begin{rem}
In \cite{davenport20141}, a similar experiment was done with movie ratings which showed a significant improvement for using 1-bit data instead of original (multi-bit) data. However, here we just see a small improvement. We believe this is due to the very careful way the in-car music data was gathered. To be precise 1-bit matrix (and tensor) completion seems to be working better when there is an implicit noise in the ratings given by the users which can't be accounted for when we fit the data exactly.
\end{rem}
\begin{rem}
It is generally true that the prediction should be more accurate when the original rating is further away from the average. This trend can still be seen in Table \ref{table:incarmusic} as well except for one case which is due to the very few instances of 0-rating in the test data (generally 3 to 5 ratings).
\end{rem}

\subsection{A Restaurant Context-Aware Recommender System}
Next, we apply our algorithm to a restaurant recommender system \cite{ramirez2014post}. The data was gathered by asking 50 users who went to 40 restaurants in the city of Tijuana and rated their experience through a Web based platform. The contextual information added is the day of the week (midweek and weekend) and the place (school, home and work). The resulting tensor is in $\mathbb{R}^{40 \times 50 \times 6}$ and the ratings are ranged from 1 to 5 with an average of 3.40. Similar to the In-car music recommender system we take the 1-bit information to be whether or not a rating is higher than the average or not. In Table \ref{table:restaurnat_percent} we report the percentage of correct above-or-below-average predictions for different ratings in the test set which we will refer to as \textit{1-bit predictions}. To be precise we show what percentage of the test ratings have been correctly predicted to be above or below average. In Table \ref{table:restaurnat_MEA} we show the mean absolute error of the recovered tensor on the test set, i.e., taking $\Omega_t$ to be the test set, we show $MAE:=\frac{\sum_{\omega \in \Omega_t}|T^{\sharp}(\omega)-\hat{T}(\omega)|}{|\Omega_t|}$. For this restaurant-ratings data set, a naive matricization, leading to a $50 \times 240$ matrix, results in an overall prediction accuracy of 60\%. A more complicated matrix completion version would have repeated user-restaurant ratings (with different context information) and  it is unclear how to compare it with tensor completion. Instead, we focus on a few important questions regarding 1-bit tensor completion.
\begin{enumerate}
\item The first question is the importance of using max-qnorm regularizer. In the first two rows of Tables \ref{table:restaurnat_percent} and \ref{table:restaurnat_MEA} we show the results without using max-qnorm and just using tensor factorization with $r=10$, i.e., by doing alternating minimization over $N_i \times 10$ factors ($r=10$ is the best rank found by numerous  empirical experiments). Both results of 1-bit prediction and the mean absolute error is generally worse than using max-qnorm. The difference is significant in the 1-bit predictions. It is worth mentioning though that the tensor factorization is a non-convex problem and a more complex algorithm might get better results but it is not the focus of this paper.\medskip
\item The second question we investigate is the effect of choosing $f$ (in \ref{negative_loglikelihood}) on the results. The 4th and 5th row of Tables \ref{table:restaurnat_percent} and \ref{table:restaurnat_MEA} shows the results for two different values of $\sigma$. In short, smaller noise (not too small though) result in better 1-bit predictions and larger noise (not too large though) results in better prediction of the actual score. This behavior can be explained by examining the likelihood function \ref{table:restaurnat_MEA} more closely. When $\sigma$ is small the dithering function $f$ is spiky around $X(\omega)=0$ and therefore it is more sensitive to whether or not $X(\omega)$ is positive or negative rather not how large it is. Therefore, using larger values of $\sigma$ does a better job in recovering the original rating but has less sensitivity to the sign of $X(\omega)$. This can be seen in Tables \ref{table:restaurnat_percent} and \ref{table:restaurnat_MEA} where we recover the sign of the ratings better when we use $\sigma=0.1$ but do worse in terms of mean absolute error. This is more evident in Table \ref{table:restaurnat_MEA} where using smaller $\sigma$ does a better job in recovering the ratings that are close to the average and worse in the ratings that are further from the average. Remark \ref{remark_sigma_value} investigates these differences in more details.\medskip
%Considering these observations, we experimented on using the following loss function which combines a small noise and a large noise in the log-likelihood function.
%\begin{equation}\label{negative_loglikelihood_2sigmas}
%\begin{aligned}
%\mathcal{L}_{\Omega,Y}(X)=&\sum_{\omega \in \Omega}\left(\mathds{1}_{[Y_\omega=1]}log(\frac{1}{f_1(X_\omega)})+\mathds{1}_{[Y_\omega=-1]}log(\frac{1}{1-f_1(X_\omega)})\right)+\\ &\sum_{\omega \in \Omega}\left(\mathds{1}_{[Y_\omega=1]}log(\frac{1}{f_2(X_\omega)})+\mathds{1}_{[Y_\omega=-1]}log(\frac{1}{1-f_2(X_\omega)})\right),
%\end{aligned}
%\end{equation}
%where we use $f_1(x)=\Phi(\frac{x}{0.2}$ as a small noise and $f_2(x)=\Phi(\frac{x}{0.2}$ as a larger noise.
%
\end{enumerate}
As expected the best results for recovering the original ratings (MAE in Table \ref{table:restaurnat_MEA}) is obtained by exact tensor completion. However, notice that 1-bit TC outperform exact tensor completion for ratings that are around the mean significantly and struggles for correct scale of the ratings that are far away from the mean.
%\begin{rem}\label{remark_sigma_value}
%Ignoring the infinity-norm constraints, in the application section, we don't dither the original tensor to obtain the 1-bit measurements and use $f(x)=\Phi(\frac{x}{\sigma})$ in the reconstruction algorithm. Therefore, in theory changing the value of $\sigma$ should not change the results of predicting the sign of the tensor. This can be observed by noticing that
%$$\Phi(\frac{T(\omega)}{\sigma_1}) = \Phi(\frac{\frac{\sigma_2}{\sigma_1} T(\omega)}{\sigma_2})$$
%which shows that the recovered tensors obtained by using the function $f(x)=\Phi(\frac{x}{\sigma_1})$ and $f(x)=\Phi(\frac{x}{\sigma_2})$ should be multiple scalars of each other and the sign predictions should be the same. Although this is not the case in our experiments, the difference in the sign predictions is small and it is caused because the algorithm might not converge to a global optimum.
%\end{rem}
%\begin{rem}
%in synthetic data we dither the data using the noise function $f$ but in applications we rely on the intrinsic noise in the user's ratings and model it with a reasonable noise function.
%\end{rem}
\begin{table}
\centering
\begin{tabular}{|r||c | c | c | c | c || c|}
\hline
Original rating & 1 & 2 & 3 & 4 & 5 & Overall \\ \hline \hline
multi-bit TC (TF) & 73\% & 58\% & 63\% & 43\% & 73\% & 74\% \\ \hline
1-bit TC (TF)  & 70\% & 67\% & 74\% & 57\% & 79\% & 77\% \\ \hline
multi-bit TC (max-qnorm)  & \bf{85\%} & \bf{72\%} & 69\% & 63\% & 89\% & 81\% \\ \hline
1-bit TC (max-qnorm, $\sigma=0.1$) & 80\% & \bf{72\%} & \bf{72\%} & \bf{69\%} & \bf{92\%} & \bf{84\%} \\ \hline
1-bit TC (max-qnorm, $\sigma=0.5$)& 79\% & 69\% & 69\% & 66\% & 91\% & 82\% \\ \hline
\end{tabular}
\caption{\small \sl Results of a comparison between 1-bit and multi-bit matrix tensor completion algorithms on Tijuana restaurant data \cite{ramirez2014post} for {\bf{predicting whether the unobserved ratings were above or below average}}. TF refers to using tensor factorization without the max-qnorm regularization. For the 1-bit results we use $f(x)=\Phi(\frac{x}{\sigma})$.}
 \label{table:restaurnat_percent}
\end{table}

\begin{table}
\centering
\begin{tabular}{|r||c | c | c | c | c || c|}
\hline
Original rating & 1 & 2 & 3 & 4 & 5 & Overall \\ \hline \hline
multi-bit TC (TF) & 1.80 & 1.22 & 0.74 & 0.88 & 1.10 & 0.97 \\ \hline
1-bit TC (TF)  & 1.90 & 1.02 & 0.62 & 0.85 & 0.96 & 0.97 \\ \hline
multi-bit TC (max-qnorm)  & \bf{1.5} & \bf{1.01} & 0.54 & 0.64 & \bf{0.68} & \bf{0.76} \\ \hline
1-bit TC (max-qnorm, $\sigma=0.1$) & 2.20 & 1.25 & \bf{0.29} & \bf{0.46} & 1.25 & 1.11 \\ \hline
1-bit TC (max-qnorm, $\sigma=0.5$)& 1.52 & 1.12 & 1.08 & 1.14 & 0.84 & 1.02 \\ \hline
\end{tabular}
\caption{\small \sl Results of a comparison between 1-bit and multi-bit matrix tensor completion algorithms on Tijuana restaurant data \cite{ramirez2014post} {\bf{showing the mean absolute error}}.}
 \label{table:restaurnat_MEA}
\end{table}
\begin{rem}
Although in theory the sign predictions should be the same with $\sigma=0.1$ and $\sigma=0.5$, we see a slight difference in the sign predictions in Table \ref{table:restaurnat_percent}. We believe this difference is a result of the algorithm not converging to the global optimum.
\end{rem}
\section{Proofs}\label{section_proofs}

\subsection{Proof of Theorem \ref{theorem_maxnorm}}
Our proof closely follows the one in \cite[Section 7.1]{cai2013max}. Therefore, we just briefly explain the steps. In what follows we prove the max-qnorm constrained ML estimation error. The proof for M-norm constraint one \eqref{upper_bound_Mnorm} follows the exact same steps where the only difference is a constant difference in the Rademacher complexity of the unit balls of these two norms.
Define the loss function $g(x;y):\mathbb{R} \times \{\pm 1\} \to \mathbb{R}$ as:
$$g(x;y)=\mathds{1}_{[y=1]}log(\frac{1}{f(x)})+\mathds{1}_{[y=-1]}log(\frac{1}{1-f(x)}).$$
Regarding the tensor completion problem as a prediction problem where we consider the tensor $T \in \mathbb{R}^{N_1 \times N_2 \times \cdots \times N_d}$ as a function from $[N_1] \times [N_2] \times \cdots [N_d] \to \mathbb{R}$ where the value of the function at $(i_1,\cdots,i_d)$ is the corresponding entry, the proof is based on a general excess risk bound developed in \cite{bartlett2002rademacher}.\newline
For a subset $S=\{\omega_1,\omega_2,\cdots,\omega_m\}$ of the  set $[N_1]\times \cdots\times [N_d]$ and a tensor $T$, $D_S(T;Y)$ is the average empirical loss according to $g$. Precisely, $D_S(T;Y)=\frac{1}{m}\sum_{i=1}^m g(T(\omega_i);Y(\omega_i))$. When the sample set $S$ is drawn i.i.d. according to $\Pi$(with replacement), we define:
$$D_{\Pi}(T;Y) := \mathbb{E}_{S \sim \Pi}[g(T(\omega_i);Y(\omega_i))]=\sum_{i=1}^m \pi_{\omega_i} g(T(\omega_i),Y(\omega_i)).$$
Notice that since $\hat{T}_{\m}$ is the optimal solution of optimization problem \eqref{optimization_maxnorm} and $T^{\sharp}$ is feasible for this problem we have:
\begin{equation}
D_S(\hat{T}_{\m};Y) \leq D_S(T^{\sharp};Y) = \frac{1}{m}\sum_{i=1}^m g(T^{\sharp}(\omega_i);Y(\omega_i))
\end{equation}
Therefore, we have
\begin{equation}\label{long_inequality}
\begin{aligned}
\mathbb{E}_{Y}&[D_{\Pi}(\hat{T}_{\m};Y)-D_{\Pi}(T^{\sharp};Y)]\\
&=\mathbb{E}_{Y}[D_{\Pi}(\hat{T}_{\m};Y)-D_{S}(T^{\sharp};Y)]+\mathbb{E}_{Y}[D_{S}(T^{\sharp};Y)-D_{\Pi}(T^{\sharp};Y)]\\
&\leq \mathbb{E}_{Y}[D_{\Pi}(\hat{T}_{\m};Y)-D_{S}(\hat{T}_{\m};Y)]+\mathbb{E}_{Y}[D_{S}(T^{\sharp};Y)-D_{\Pi}(T^{\sharp};Y)]\\
&\leq \underset{T \in K^T_{\m}(\alpha,R_{\m})} {\text{sup}} \lbrace \mathbb{E}_{Y}[D_{\Pi}(T;Y)]-\mathbb{E}_{Y}[D_{S}(T;Y)] \rbrace+\mathbb{E}_{Y}[D_{S}(T^{\sharp};Y)-D_{\Pi}(T^{\sharp};Y)].
\end{aligned}
\end{equation}
Notice that the left hand side of \eqref{long_inequality} is equivalent to weighted Kullback-Leibler divergence between $f(T^{\sharp})$ and $f(\hat{T}_{\m})$.\newline
Now we focus on bounding the right hand side of \eqref{long_inequality}. Using Hoeffding's inequality on the random variable $Z_\omega := g(T^{\sharp}(\omega);Y(\omega))-\mathbb{E}[g(T^{\sharp}(\omega);Y(\omega))]$ we conclude that with probability $1-\delta$ over choosing the sampling subset $S$:
\begin{equation}\label{risk_for_actualT*}
D_{S}(T^{\sharp};Y)-D_{\Pi}(T^{\sharp};Y) \leq U_{\alpha} \sqrt{\frac{\log(\frac{1}{\delta})}{2m}}.
\end{equation}
Moreover, a combination of Theorem 8, (4) of Theorem 12 from \cite{bartlett2002rademacher} and the upper bound \eqref{Rademacher_maxnorm} and noting that $g(\dot,y)$ is an $L_{\alpha}$-Lipschitz function yields:
\begin{equation}\label{risk_sup_inK}
\underset{T \in K^T_{\m}(\alpha,R_{\m})} {\text{sup}} \lbrace \mathbb{E}_{Y}[D_{\Pi}(T;Y)]-\mathbb{E}_{Y}[D_{S}(T;Y)] \rbrace \leq 12 L_{\alpha} R_{\m} c_1 c_2^d \sqrt{\frac{dN}{m}}+U_{\alpha} \sqrt{\frac{8\log(\frac{2}{\delta})}{m}}.
\end{equation}
Gathering \eqref{long_inequality}, \refeq{risk_for_actualT*}, and \eqref{risk_sup_inK} we get:
\begin{equation}\label{kullback_bound}
\mathbb{K}_{\Pi}(f(T^{\sharp})||f(\hat{T}_{\m})) \leq 12 L_{\alpha} R_{\m} c_1 c_2^d \sqrt{\frac{dN}{m}}+U_{\alpha} \sqrt{\frac{8\log(\frac{2}{\delta})}{m}} + U_{\alpha} \sqrt{\frac{\log(\frac{1}{\delta})}{2m}}.
\end{equation}
This, together with \eqref{KL_div} and \cite[Lemma 2]{davenport20141} proves \eqref{upper_bound}. The upperbound \eqref{upper_bound_Mnorm} can be proved by following the exact same arguments with the only difference being in the right hand side of \eqref{risk_sup_inK} where the Rademacher complexity of unit M-norm ball does not have the constant $c_1 c_2^d$ and the upper bound $R_M$ is different than $R_{\m}$.
\subsection{Proof of Theorem \ref{theorem_information_theoritic}}\label{section_proof_theorem_upperbound}
\begin{lem}[{\cite[Lemma 18]{ghadermarzy2017near}}]\label{packing}
Let $r=\floor{(\frac{R}{\alpha K_G})^2}$ and let $K_{M}^T(\alpha,R)$ be defined as in Section \ref{section_information} and let $\gamma \leq 1$ be such that $\frac{r}{\gamma^2}$ is an integer and suppose $\frac{r}{\gamma^2} \leq N$. Then the there exist a set $\chi_{(\alpha,\gamma)}^T \subset K_{M}^T(\alpha,R) $ with
$$|\chi_{(\alpha,\gamma)}^T| \geq \exp{\frac{rN}{16\gamma^2}}$$
such that
\begin{enumerate}
\item For $T \in \chi_{(\alpha,\gamma)}^T$, $|T(\omega)|=\alpha \gamma$ for $\omega \in \{[N] \times [N] \cdots [N]\}$.
\item For any $T^{(i)},\ T^{(j)}\ \in \chi_{(\alpha,\gamma)}^T$, $T^{(i)} \neq T^{(j)}$
$$\|T^{(i)}-T^{(j)}\|_F^2 \geq \frac{\alpha^2 \gamma^2 N^d}{2}$$
\end{enumerate}
\end{lem}
Next choosing $\gamma$ in a way that $\frac{r}{\gamma^2}$ is an integer and
$$ \frac{4\sqrt{2} \epsilon}{\alpha} \leq \gamma \leq \frac{8\epsilon}{\alpha}$$
we transform $\chi_{(\alpha,\gamma)}^T$ so that the entries of each element of the packing set come from the set $\{\alpha,\alpha^{'}:=(1-\gamma)\alpha)\}$ by defining
$$\chi := \{X + \alpha (1-\frac{\gamma}{2}) {\bf{1}} : X \in \chi_{(\frac{\alpha}{2},\gamma)}^T\}$$
Notice that for any $X \in \chi$, $\|X\|_{\infty} \leq \alpha$ and $\|X\|_{M} \leq R_M$ and therefore, we can use $\chi$ as a packing set for $K_M^T$. Once the packing sets are generated, using the exact same proof as \cite[Section A.3]{davenport20141}, we can prove an upper bound. Choose $\epsilon$ to be
\begin{equation}\label{epsilon_upperbound}
\epsilon^2 = \text{min} \left(\frac{1}{1024}, C_2 \alpha \sqrt{\beta_{\frac{3\alpha}{4}} }\sqrt{\frac{rN}{m}} \right).
\end{equation}
And assume there is an algorithm that using 1-bit measurements $Y_{\Omega}$ of $X$, returns $\hat{X}$ such that
$$\frac{1}{N^d} \|X-\hat{X} \|_F^2 \leq \epsilon^2,$$
with probability at least $\frac{1}{4}$. Notice that this means that with probability at least $\frac{1}{4}$, we can distinguish the elements of the packing set as for any $X^{(i)},X^{(j)} \in \chi$, $\frac{1}{N^d} \|X^{(i)} - X^{(j)} \|_F^2 \geq 4 \epsilon^2$. Hence, defining $X^{\ast}$ to be the closest element of $\chi$ to $\hat{X}$, if we show that $\mathbb{P}(X \neq X^{\ast}) > \frac{3}{4}$ the proof is done. Otherwise, using Fano's inequality \cite{davenport20141} and following the steps in Section A.3.2 of \cite{davenport20141} we have
$$\frac{1}{4} \leq 1- \mathbb{P}(X \neq X^{\ast}) \leq 16 \gamma^2 \left(\frac{\frac{64m\epsilon^2}{\beta_{\alpha^{\prime}}}+1}{rN}\right) \leq 1024 \epsilon^2 \left(\frac{\frac{64m\epsilon^2}{\beta_{\alpha^{\prime}}}+1}{\alpha^2 rN}\right),$$
where $\alpha^{\prime}:=(1-\gamma) \alpha$.
Comparing $64m\epsilon^2$ and $\beta_{\alpha^{\prime}}$, either

$$\frac{1}{4} < \frac{2048 \epsilon^2}{\alpha^2 r N},$$
which is a contradiction when $C_0 > 8$ or
$$\epsilon^2 > \frac{\alpha\sqrt{\beta_{\alpha^{\prime}}}}{512 \sqrt{2}} \sqrt{\frac{rN}{m}},$$
which is a contradiction when $C_2 \leq \frac{1}{512 \sqrt{2}}$.
\bibliography{sparse}
\bibliographystyle{plain}
\begin{appendices}
\section{Rademacher Complexity}\label{section_appendix_radamache}
\begin{dfn}\cite{cai2013max}
Let $\mathbb{P}$ be a probability distribution on a set $\chi$ and assume the set $S:=\lbrace X_1,\cdots, X_m\rbrace$ is $m$ independent samples drawn from $\chi$ according to $\mathbb{P}$. For a class of Functions $\mathbb{F}$ defined from $\chi$ to $\mathbb{R}$, its empirical Rademacher complexity is defined as:
$$\hat{R}_S(\mathbb{F})=\frac{2}{|S|}\mathbb{E}_{\epsilon}[\underset{f \in \mathbb{F}}{sup}|\sum_{i=1}^{i=m}\epsilon_if(X_i)|],$$
where $\epsilon_i$ is a Radamacher random variable. Moreover, the Rademacher complexity with respect to distribution $\mathbb{P}$ over a sample $S$ of $|S|$ point drawn independently from $S$ is defined as the expectation of the empirical Rademacher complexity defined as:
$$R_{|S|}(\mathbb{F})=\mathbb{E}_{S\sim \mathbb{P}}[\hat{R}_S(\mathbb{F})]$$
\end{dfn}
The following two lemmas were proved in \cite{ghadermarzy2017near}.
\begin{lem}\label{Rademacher_atomicnorm}
$\underset{S:|S|=m}{sup} \hat{R}_S(\mathbb{B}_{M}(1)) < 6 \sqrt{\frac{dN}{m}}$ 
\end{lem}

\begin{lem}\label{Rademacher_maxnorm}
$\underset{S:|S|=m}{sup} \hat{R}_S(\mathbb{B}_{\m}^T(1)) < 6 c_1 c_2^d \sqrt{\frac{dN}{m}}$ 
\end{lem}

\section{Discrepancy}\label{section_appendix_discrepency}
In this we briefly describe some more notation which are direct generalization of corresponding definitions for matrix distances. First, the Hellinger distance for two scalars $0 \leq p,q \leq 1$ is defined by:
$$d_H^2(p,q)=(\sqrt{p}-\sqrt{q})^2+((\sqrt{1-p}-\sqrt{1-q})^2,$$
which gives a standard notation on the distance between two probability distributions. We generalize this definition to define Hellinger between two tensors $T,U \in [0,1]^{N_1 \times \cdot N_d}$ as:
$$d_H^2(T,U)=\frac{1}{N_1 \cdots N_d}\sum_{i_1,\cdots,i_d}d_H^2(S(i_1,\cdots,i_d),U(i_1,\cdots,i_d))$$
The Kullback-Leibler divergence between two tensors $T,U \in [0,1]^{N_1 \times \cdot N_d}$ is generalized as:
$$\mathbb{K}(T||U):=\frac{1}{N_1 \cdots N_d}\sum_{i_1,\cdots,i_d}\mathbb{K}(S(i_1,\cdots,i_d)||U(i_1,\cdots,i_d)),$$
where $\mathbb{K}(p||q):=p \log(\frac{p}{q}) + (1-p) \log (\frac{1-p}{1-q})$, for $0 \leq p,q \leq 1$. It is easy to prove that for any two scalars $0 \leq p,q \leq 1$, $d_H^2(p,q) \leq \mathbb{K}(p||q)$ and hence:
\begin{equation}\label{KL_div}
d_H^2(T,U) \leq \mathbb{K}(T||U).
\end{equation}

\end{appendices}
\end{document}